\documentclass[10pt]{amsart}
\usepackage{amscd}
\usepackage{amssymb}

\newtheorem{thm}[equation]{Theorem}
\newtheorem{cor}[equation]{Corollary}
\newtheorem{prop}[equation]{Proposition}
\newtheorem{lem}[equation]{Lemma}

\theoremstyle{definition}
\newtheorem{dfn}[equation]{Definition}
\newtheorem{rem}[equation]{Remark}
\newtheorem{exa}[equation]{Example}

\newtheorem{que}[equation]{Question}
\numberwithin{equation}{section}

\newcommand{\opn}{\operatorname}

\newcommand{\cat}[1]{\operatorname{\mathsf{#1}}}

\newcommand{\rmitem}[1]{\item[\text{\textup{(#1)}}]}

\newcommand{\mfrak}[1]{\mathfrak{#1}}

\newcommand{\msf}[1]{\mathsf{#1}}

\newcommand{\mrm}[1]{\mathrm{#1}}
\newcommand{\mbb}[1]{\mathbb{#1}}

\newcommand{\tup}[1]{\textup{#1}}

\newcommand{\Ext}{\opn{Ext}}

\newcommand{\fp}{{\mathfrak p}}

\title[Dualizing Complexes and Tilting Complexes]
{Dualizing Complexes and Tilting Complexes over Simple Rings}
\author{Amnon Yekutieli and James J.\ Zhang }
\address{Department of Mathematics,
Ben Gurion University, Be'er Sheva 84105, ISRAEL}
\email{amyekut@math.bgu.ac.il}
\address{Department of Mathematics, University of
Washington, Seattle, WA 98195, USA}
\email{zhang@math.washington.edu}
\date{27 September 2001; v1/d4}

\begin{document}
\maketitle


\setcounter{section}{-1}
\section{Introduction}

Simple rings, like fields, are literally `simple' in many ways. 
Hence quite a few invariants of rings become trivial for simple rings. 
We show that this principle applies to the derived Picard group, 
which classifies dualizing complexes over a ring.

In this paper all rings are algebras over a base field $k$, ring 
homomorphisms are all over $k$, and bimodules are all $k$-central. 
The symbol $\otimes$ denotes $\otimes_{k}$. For a ring $B$, 
$B^{\circ}$ denotes the opposite ring. 

We shall write $\cat{Mod} A$ for the category of left 
$A$-modules, and $\msf{D}^{\mrm{b}}(\cat{Mod} A)$ will stand for 
the bounded derived category. A brief review of key 
definitions such as dualizing complexes, two-sided tilting 
complexes and the derived Picard group $\opn{DPic}(A)$ is included in 
the body of the paper.

\begin{thm}
\label{yy0.1}
Let $A$ and $B$ be rings and let 
$T\in \msf{D}^{\mrm{b}}(\cat{Mod}(A\otimes B^\circ))$ 
be a two-sided tilting complex. Suppose either $A$ or $B$ is a 
Goldie simple ring. 
\begin{enumerate}
\item{} $T\cong P[n]$ for some integer $n$ and 
some invertible $A$-$B$-bimodule $P$. Therefore $A$ and $B$ are 
Morita equivalent, and in particular both are Goldie simple rings.
\item{} The structure of the derived Picard group of $A$ is
$\opn{DPic}(A) = {\mathbb Z}\times \opn{Pic}(A)$.
\end{enumerate}
\end{thm}

An algebra is called {\it Gorenstein} if it has finite left
and right injective dimension. 

\begin{thm}
\label{yy0.2}
Let $A$ be a left noetherian ring and let $B$ be a right noetherian
ring. Let $R$ be a dualizing complex over $(A,B)$.
Assume either of the two conditions below hold.
\begin{enumerate}
\rmitem{i} $A$ and $B$ are both Goldie simple rings.
\rmitem{ii} Either $A$ or $B$ is a Goldie simple ring, and 
either $A$ or $B$ is noetherian and admits some dualizing 
complex.
\end{enumerate}
Then $R \cong P[n]$ for some integer $n$ and some 
invertible $A$-$B$-bimodule $P$, the rings $A$ and $B$ are Morita 
equivalent, and both are noetherian Gorenstein simple rings.
\end{thm}

One motivating question is to classify all dualizing and tilting 
complexes over the Weyl algebras. When the base field has characteristic
zero, this question is answered by Theorems \ref{yy0.1} and \ref{yy0.2}.
When the base field has positive characteristic, the same answer is given 
in Section 5.

Theorem \ref{yy0.2} also has a surprising consequence.

\begin{cor}
\label{yy0.3} Let $A$ be a filtered ring such that the associated 
graded ring $\opn{gr} A$ is connected graded and noetherian. Suppose
either one of the following conditions holds:
\begin{enumerate}
\rmitem{i} $\opn{gr} A$ is commutative.
\rmitem{ii} $\opn{gr} A$ is PI.
\rmitem{iii} $\opn{gr} A$ is FBN. 
\rmitem{iv} $\opn{gr} A$ has enough
normal elements in the sense of \cite[p.36]{YZ}.
\rmitem{v} $\opn{gr} A$ is a factor ring of a graded AS-Gorenstein 
ring. 
\end{enumerate}
If $A$ is simple, then $A$ is Gorenstein. In cases \tup{(i-iv)}, 
$A$ is also Auslander-Gorenstein and Cohen-Macaulay.
\end{cor}

For example, every simple factor ring $A$ of the enveloping
algebra $U(L)$ of a finite dimensional Lie algebra $L$
is Auslander-Gorenstein and Cohen-Macaulay. 
This is also true for simple factor rings of many
quantum algebras listed in \cite{GL}.

In Section 1 we review some basic facts about bimodules over simple 
rings. Theorem \ref{yy0.1} is proved at the end of Section 2.
Theorem \ref{yy0.2} is proved in Section 3, and Corollary 
\ref{yy0.3} is proved in Section 4. In 
Section 5 we prove statements analogous to Theorems \ref{yy0.1} and
\ref{yy0.2} when $A$ is a Weyl algebra over a base field $k$ of 
positive characteristic. In Section 6 we discuss an example of Goodearl 
and Warfield which shows that not every noetherian simple ring is 
Gorenstein.

\section{Preliminaries}

Let $A$ be a ring (i.e.\ a $k$-algebra). 
By an $A$-module we mean a left $A$-module. 
With this convention an $A^\circ$-module means a right $A$-module.
A finitely generated $A$-module is called {\it finite}. 

Our reference for derived categories is \cite{RD}. As for 
derived categories and derived functors of bimodules, such as 
$\mrm{R} \opn{Hom}$ and $\otimes^{\mrm{L}}$,
the reader is referred to \cite{Ye1} and \cite{Ye2}.

The following elementary facts will be used later.

\begin{lem}
\label{yy1.1}
Let $A$ be a ring and let $B$ be a \tup{(}left and right\tup{)} 
Goldie simple ring. Let $M$ be a nonzero $A$-$B$-bimodule 
finite on both sides. Then:
\begin{enumerate}
\item{} $M$ is a generator of $\cat{Mod} B^{\circ}$.
\item{} If the canonical homomorphism
$A \to \opn{End}_{B^{\circ}}(M)$ 
is bijective, then $M$ is projective as $A$-module.
\item{} Suppose that $A$ is also a Goldie simple ring, and that 
both $A \to \opn{End}_{B^{\circ}}(M)$ 
and
$B^{\circ} \to \opn{End}_{A}(M)$ 
are bijective. Then $M$ is an invertible bimodule.
\end{enumerate}
\end{lem}

\begin{proof} (1) 
Suppose 
$M = \sum_{i = 1}^{p} A \cdot m_{i}$
and let $N_{i} := \opn{Ann}_{B^{\mrm{o}}} (m_{i})$.
Then
$\opn{Ann}_{B^{\mrm{o}}} M = \bigcap_{i = 1}^{p} N_{i}$.
Since $B$ is a simple ring and $M \neq 0$ we must have 
$\opn{Ann}_{B^{\mrm{o}}} M = 0$. 
Hence for some $i$ the right ideal $N_{i} \subset B$ is not 
essential. This implies the element $m_{i}$ is not torsion, and so 
the $B^{\mrm{o}}$-module $M$ is not torsion.

At this point we can forget the $A$-module structure on $M$. So 
let $M$ be a finite $B^{\mrm{o}}$-module that is not torsion. 
We will show that
$\opn{Hom}_{B^{\mrm{o}}}(M, B) \neq 0$.
Replacing $M$ by a quotient of it we may assume $M$ is a finite 
uniform torsion-free $B^{\mrm{o}}$-module. In this case we have 
injections
$M \to M \otimes_{B} Q \to Q$
where $Q$ is the total ring of fractions of $B$.

Without loss of generality we can assume $M$ is a finite
$B^{\mrm{o}}$-submodule of $Q$. Thus 
$M = \sum_{i = 1}^{q} s_{i}^{-1} x_{i} \cdot B^{\mrm{o}}$
for certain $s_{i}, x_{i} \in B$ with $s_{i}$ regular elements. 
Passing to a left common denominator we have
$s_{i}^{-1} x_{i} = s^{-1} y_{i}$
for suitable $s, y_{i} \in B$. Therefore left multiplication by $s$
is a nonzero $B^{\mrm{o}}$-linear map
$\lambda_{s}: M \to B$. 

Finally we reduce to the case of a finite $B^{\mrm{o}}$-module $M$ 
such that $\opn{Hom}_{B^{\mrm{o}}}(M, B) \neq 0$.
Let $I \subset B$ be the union of the images of all 
$B^{\mrm{o}}$-linear homomorphisms $M \to B$. This is a nonzero 
two-sided ideal, and hence $I = B$. So there are some 
homomorphisms
$\phi_{i} : M \to B$ such that 
$1 \in \sum_{i = 1}^{r} \phi_{i}(M) \subset B$. Thus
$\sum \phi_{i} : M^{r} \to B$ is surjective, proving that $M$ is a 
generator of $\cat{Mod} B^{\mrm{o}}$.

\medskip \noindent
(2) Since $M$ is a generator of $\cat{Mod} B^{\circ}$
and $A \to \opn{End}_{B^{\circ}}(M)$ is bijective,
a theorem of Morita \cite{Mo} (see \cite[17.8]{AF}) says that 
$M$ is a finite projective $A$-module.

\medskip \noindent
(3) By parts (1) and (2), the $A$-module $M$ is also a finite
projective generator. By Morita's theorem the bimodule $M$ 
is invertible. 
\end{proof}

\begin{lem}
\label{yy1.2} 
Let $M$ be a bounded complex of $B^\circ$-modules with nonzero 
cohomology such that $\Ext^i_{B^\circ}(M,M)=0$ for all $i< 0$. 
Let $i_0 := \min\{i\;|\; \mrm{H}^i M\neq 0\}$ and 
$j_0 := \max\{j\;|\; \mrm{H}^j M \neq 0\}$.
If $i_0\neq j_0$ \tup{(}i.e.\ $i_0<j_0$\tup{)}, 
then 
$\opn{Hom}_{B^\circ} (\mrm{H}^{j_0} M,\mrm{H}^{i_0} M)=0$.
\end{lem}

\begin{proof} This is true because a nonzero morphism from 
$\mrm{H}^{j_0} M$ to $\mrm{H}^{i_0} M$ gives rise to a nonzero element
in $\Ext^{i_0-j_0}_{B^\circ}(M,M)$.
\end{proof}

\begin{lem} 
\label{yy1.3}
Let $M$ be a bounded complex of $A$-$B$-bimodules with nonzero 
cohomology. Suppose the following conditions hold:
\begin{enumerate}
\rmitem{i} $B$ is Goldie and simple.
\rmitem{ii} $\Ext^i_{B^\circ}(M,M)=0$ for all $i\neq 0$.
\rmitem{iii} $\mrm{H}^{j_0} M$ is finite on both sides, where 
$j_0$ is as in Lemma \tup{\ref{yy1.2}}.
\end{enumerate}
Then $M\cong (\mrm{H}^{j_0} M)[-j_0]$ in
$\msf{D}(\cat{Mod} (A \otimes B^{\circ}))$.
\end{lem}

\begin{proof} By Lemma \ref{yy1.1}(1), $\mrm{H}^{j_0} M$ is a generator
of $\cat{Mod} B^{\circ}$. Let $i_0$ be as in Lemma \ref{yy1.2}.
If $i_0<j_0$ then the conclusion of Lemma \ref{yy1.2}
contradicts the fact that $\mrm{H}^{j_0} M$ is a generator of
$\cat{Mod} B^{\circ}$. Therefore $i_0=j_0$ and the assertion follows. 
\end{proof}


\section{Two-Sided Tilting Complexes}

The following definition is due to Rickard \cite{Ri1,Ri2} and 
Keller \cite{Ke}. Recall that ``ring'' means ``$k$-algebra''.

\begin{dfn}
\label{yy2.1}
Let $A$ and $B$ be rings and let 
$T \in  \msf{D}^{\mrm{b}}(\cat{Mod}(A \otimes B^{\circ}))$ 
be a complex. We say $T$ is a {\it two-sided tilting complex} over 
$(A, B)$ if there exists a complex 
$T^{\vee} \in  \msf{D}^{\mrm{b}}(\cat{Mod}(B \otimes A^{\circ}))$
such that 
$T \otimes^{\mrm{L}}_{B} T^{\vee} \cong A$
in $\msf{D}(\cat{Mod}(A \otimes A^{\circ}))$, 
and $T^{\vee} \otimes^{\mrm{L}}_{A} T \cong B$
in $\msf{D}(\cat{Mod}(B \otimes B^{\circ}))$.
\end{dfn}

The complex $T$, when considered as a complex of left $A$-modules, 
is perfect, and the set 
$\cat{add} T \subset \msf{D}^{\mrm{b}}(\cat{Mod} A)$, namely
the direct summands of finite direct sums of $T$, generates the 
category $\msf{D}^{\mrm{b}}(\cat{Mod} A)_{\mrm{perf}}$ 
of perfect complexes. The formula for $T^{\vee}$ is
$T^{\vee} \cong \mrm{R} \opn{Hom}_{A}(T, A)$. 
The adjunction morphism
$B \mapsto \mrm{R} \opn{Hom}_{A}(T, T)$
in $\msf{D}(\cat{Mod}(B \otimes B^{\circ}))$
is an isomorphism. The functor
$M \mapsto T \otimes^{\mrm{L}}_{B} M$ is an equivalence
$\msf{D}(\cat{Mod} B) \to \msf{D}(\cat{Mod} A)$
preserving boundedness. By symmetry there are three more 
variations of all these assertions (e.g.\ $T^{\vee}$ is a perfect 
complex of $A^{\circ}$-modules). See \cite{Ye2} for proofs. 

The next definition is due to the first author \cite{Ye2}. When 
$B = A$ we write $A^{\mrm{e}} := A \otimes A^{\mrm{\circ}}$.

\begin{dfn}
\label{yy2.2} 
Let $A$ be ring. The {\it derived Picard group} of $A$
is defined to be
\[ \opn{DPic}(A) := \frac
{\{\text{ two-sided tilting complexes }
T \in \msf{D}^{\mrm{b}}(\cat{Mod} A^{\mrm{e}}) \}}
{\text{isomorphism}} , \]
with operation $(T, S) \mapsto T \otimes^{\mrm{L}}_{A} S$.
\end{dfn}

Clearly the definition of the group $\opn{DPic}(A)$ is relative to 
the base field $k$. For instance if $A = K$ is a field extension 
of $k$ then 
$\opn{DPic}(K) = \mbb{Z} \times \opn{Gal}(K / k)$, 
where $\opn{Gal}(K / k)$ is the Galois group (cf.\ 
\cite[3.4]{Ye2}). 

The derived Picard group was computed in various cases, see 
\cite{Ye2} and \cite{MY}. As shown in \cite{Ye2}, the derived 
Picard group classifies the isomorphism classes of dualizing 
complexes (cf.\ next section). 

There are some obvious tilting complexes. If $P$ is an invertible 
$A$-bimodule and $n$ is an integer, then $T := P[n]$ is a two-sided
tilting complex. Recall that the (noncommutative) Picard group 
$\opn{Pic}(A)$ of $A$ is the group of isomorphism classes of 
invertible bimodules. It follows that $\opn{DPic}(A)$ 
contains a subgroup ${\mathbb Z}\times \opn{Pic}(A)$.

\begin{proof}[Proof of Theorem \tup{\ref{yy0.1}}] 
(1) Assume that $B$ is simple and Goldie. Let 
$$j_0 := \max\{i\;|\; \mrm{H}^i(T)\neq 0\}.$$ 
Without loss of generality 
we may assume that $j_0=0$ (after a complex shift). As in 
\cite[1.1]{Ye2}, $\mrm{H}^0(T)$ is finite on both 
sides. By Lemma \ref{yy1.3} it follows that $T\cong P$ where 
$P := \mrm{H}^0(T)$. 

Since $P$ is a two-sided tilting complex we have
\[ \opn{End}_{B^{\circ}}(P) \cong 
\mrm{H}^{0} \mrm{R} \opn{Hom}_{B^{\circ}}(T, T) \cong A . \]
By Lemma \ref{yy1.1}(2), $P$ is a projective $A$-module. According 
to \cite[2.2]{Ye2}, $P$ is an invertible 
$A$-$B$-bimodule. The functor $M \mapsto P \otimes_{B} M$ is then an 
equivalence $\cat{Mod} B \to \cat{Mod} A$.

\medskip \noindent
(2) Take $A=B$. By part (1) every tilting complex is isomorphic 
to $P[n]$. The assertion follows.
\end{proof}


\section{Dualizing Complexes}

The definition of a dualizing complex over a noncommutative graded
ring is due to the first author \cite{Ye1}. The following more 
general definition appeared in \cite{YZ}. 

\begin{dfn}
\label{yy3.1}
Assume $A$ is a left noetherian ring and $B$ is a right 
noetherian ring. A complex 
$R \in \msf{D}^{\mrm{b}}(\cat{Mod}(A\otimes B^{\circ}))$ 
is called a {\it dualizing complex over $(A,B)$} if it satisfies the 
following conditions:
\begin{enumerate}
\rmitem{i} $R$ has finite injective dimension over $A$ and over
$B^\circ$.
\rmitem{ii} $R$ has finite cohomology 
modules $A$ and over $B^\circ$. 
\rmitem{iii} The canonical morphisms $B\to \mrm{R}\opn{Hom}_A(R,R)$ in
$\msf{D}(\cat{Mod}(B\otimes B^\circ))$ and 
$A\to \mrm{R} \opn{Hom}_{B^\circ}
(R,R)$ in $\msf{D}(\cat{Mod}(A\otimes A^\circ))$ 
are both isomorphisms.
\end{enumerate}
If moreover $A=B$, we say $R$ is a {\it dualizing complex over $A$}.
\end{dfn}

Whenever we say $R$ is a dualizing complex over $(A,B)$ we are 
tacitly assuming that $A$ is left noetherian and $B$ is right 
noetherian.  

A noetherian ring $A$ is Gorenstein if and only if the bimodule $R := A$ is 
a dualizing complex. Existence of dualizing complexes
for non-Gorenstein rings is studied in \cite{VdB, YZ}. 

If $A$ is noetherian and has at least one dualizing complex then 
the derived Picard group $\opn{DPic}(A)$ classifies the 
isomorphism classes of dualizing complexes. Indeed, given a 
dualizing complex $R$, any other dualizing complex $R'$ is 
isomorphic to $R \otimes^{\mrm{L}}_{A} T$ 
for some two-sided tilting complex $T$, and $T$ is unique up to 
isomorphism.

\begin{proof}[Proof of Theorem \tup{\ref{yy0.2}}]
By Lemma \ref{yy1.3}, $R\cong P[n]$ for some bimodule $P$ and some
integer $n$. 

Since $R$ is dualizing the canonical homomorphisms
$A \to \opn{End}_{B^\circ}(P)$ and 
$B^\circ \to \opn{End}_{A}(P)$
are isomorphisms. When both $A$ and $B$ are 
Goldie and simple (condition (i)), 
Lemma \ref{yy1.1}(3) implies that $P$ is invertible.  

Now assume $A$ is noetherian, and it has some dualizing complex $R_1$
(condition (ii)).
Then by the proof of \cite[4.5]{Ye2} -- suitably modified to fit our 
situation -- the complex
$T := \mrm{R} \opn{Hom}_A(R_1,R)\in \msf{D}^{\mrm{b}}
(\cat{Mod}(A\otimes B^\circ))$ 
is a two-sided tilting complex. Since either $A$ or $B$ is a Goldie 
simple ring, it follows from Theorem \ref{yy0.1} that both $A$ and 
$B$  are Goldie simple rings. As above we deduce that $P$
is an invertible bimodule.

Under both conditions the rings $A$ and $B$ are Morita equivalent. 
Since the bimodule $P$ is a dualizing complex over $(A, B)$ 
it has finite injective 
dimension on both sides. But on the other hand $P$ is a  
progenerator on both sides, and hence $A$ has finite injective dimension 
on  the left, and $B$ has finite injective dimension on the right. 
By Morita equivalence, both $A$ and $B$ are (two-sided) noetherian 
and have finite left and right injective dimensions. 
\end{proof}

\begin{rem} 
\label{yy3.2}
One can define dualizing complexes in a slightly more general 
situation, by replacing the noetherian condition with the weaker 
coherence condition (see \cite[3.3]{Ye1}). Thus in 
Definition \ref{yy3.1} $A$ is a left coherent ring, $B$ is a right 
coherent ring, and in condition (ii) of  the word `finite' is 
replaced with `coherent'. It is not hard to check that 
a ``coherent'' version of Theorem \ref{yy0.2} holds. 
\end{rem}

\begin{exa}
\label{yy3.3}
Let $A := \bigcup A_{n}$ where $A_{n}$ is the $n$th Weyl algebra 
over the field $k$. Using the method of faithful flatness we see 
that $A$ has the following properties:
\begin{enumerate}
\rmitem{i} $A$ is neither left nor right noetherian.
\rmitem{ii} $A$ has infinite Krull, Gelfand-Kirillov, injective, and 
global dimensions.
\rmitem{iii} $A$ is a Goldie domain (i.e.\ a left and right Ore 
domain). 
\rmitem{iv} $A$ is a coherent ring.
\end{enumerate}

Suppose now $\opn{char}k =0$. Then $A$ is a simple ring. 
Therefore Theorem \ref{yy0.1} holds. For instance, the derived Picard
group of $A$ is ${\mathbb Z}\times \opn{Pic}(A)$. By the 
``coherent'' version of Theorem \ref{yy0.2} (see Remark \ref{yy3.2}) 
$A$ does not admit a dualizing complex, because $A$ is not Gorenstein. 
\end{exa}

Examples of noetherian simple rings with infinite Krull dimension
were given by Shamsuddin \cite{Sh} and Goodearl-Warfield \cite{GW2}.
It is not hard to show that these simple rings also have infinite 
injective dimension (see Section 6).



\section{The Auslander Condition}

Let $R$ be a dualizing complex over $(A,B)$ and let $M$ be 
an $A$-module. The {\it grade} of $M$ with respect to $R$ is defined to be
$$j_{R}(M)=\min\{q\; |\; \opn{Ext}^q_A(M,R)\neq 0\}.$$
The grade of a $B^\circ$-module is defined similarly.

We recall the definitions of the Auslander condition and the 
Cohen-Macaulay condition. 
Gelfand-Kirillov dimension is denoted by $\opn{GKdim}$.

\begin{dfn} \label{yy4.1}
Let $R$ be a dualizing complex over $(A,B)$.
\begin{enumerate}
\item $R$ is called {\it Auslander} if the two conditions below 
hold. 
\begin{enumerate}
\rmitem{i} For every finite $A$-module $M$, every 
$q$, and every $B^\circ$-submodule 
$N\subset \opn{Ext}^q_A(M,R)$ one has $j_{R}(N)\geq q$.
\rmitem{ii} The same holds after exchanging $A$ and $B^\circ$.
\end{enumerate}
\item If there is a constant $s$ such that 
$$j_{R}(M)+\opn{GKdim} M=s$$
for all finite $A$-modules or finite $B^\circ$-modules $M$, then $R$ 
is  called {\it Cohen-Macaulay}.
\end{enumerate}
\end{dfn}

The {\it canonical dimension} with respect to an Auslander
dualizing complex $R$ is defined to be
$$\opn{Cdim}_{R} M=-j_{R}(M)$$
for all finite $A$-modules or $B^\circ$-modules $M$.
By \cite[2.10]{YZ}, $\opn{Cdim}_{R}$ is a finitely partitive, 
exact dimension function. See \cite[6.8.4]{MR} for the definition
of dimension function. 

When $A$ is a Gorenstein ring and the bimodule $R := A$ is an
Auslander dualizing complex then $A$ is called an {\em 
Auslander-Gorenstein} ring. If $A$ is an Auslander-Gorenstein ring such 
that  $R := A$ is also Cohen-Macaulay, then $A$ is called an 
{\it Auslander-Gorenstein Cohen-Macaulay ring}.
This is the usage in \cite{Le} and \cite{Bj}. 
We remind the critical reader that unlike commutative rings, a 
noncommutative Gorenstein ring need not be either Auslander or Cohen-Macaulay.

The following is easy.

\begin{lem}
\label{yy4.3}
Let $A$, $B$, $C$ be rings. Let $L, M, N$ be bounded complexes
over $B^\circ, A, A \otimes B^\circ$ respectively,
and let $P$ be a invertible $B$-$C$-bimodule.  
\begin{enumerate}
\item{}
For every $i$ there is an isomorphism of $A$-modules
$$\Ext^i_{B^\circ}(L,N)\cong \Ext^i_{C^\circ}(L\otimes_B P,N\otimes_B P).$$
\item{}
Suppose $A$ is left noetherian and $H^jM$ is finite over $A$ for all $j$. 
Then, for every $i$ there 
is an isomorphism of $C^\circ$-modules
$$\Ext^i_{A}(M,N\otimes_B P)\cong \Ext^i_{A}(M,N)\otimes_B P.$$
\end{enumerate}
\end{lem}

\begin{proof} (1) This follows from the fact that $-\otimes_B P$ 
induces a Morita equivalence.

(2) This is obvious when $M=A[i]$. Then the assertion follows from 
the facts that $M$ has a bounded above resolution by finite free 
$A$-modules and $P$ is a flat $B^\circ$-module.
\end{proof}

\begin{prop}
\label{yy4.4} Let $A$, $B$, $C$ be rings. Let $P$ be an invertible
$B$-$C$-bimodule and $n$ an integer. 
Suppose $R$ is a dualizing complex over $(A,B)$, 
and let $R_1=R\otimes_B P[n]$, which is a dualizing complex
over $(A, C)$. Then $R$ is Auslander
\tup{(}resp.\ Cohen-Macaulay\tup{)} if and only if $R_1$ is.
\end{prop}

\begin{proof} Without loss of generality we may assume $n=0$.

Let us assume $R_1$ is Auslander; we will prove that $R$ is also 
Auslander. Given a finite $A$-module
$M$ and an integer $i$, let $N$ be a $B^\circ$-submodule of 
$\Ext^i_A(M,R)$. Then $N\otimes_B P$ is a $C^\circ$-submodule of 
$$\Ext^i_A(M,R)\otimes_A P\cong\Ext^i_A(M,R\otimes_B P)=\Ext^i_A(M,R_1).$$
By the Auslander condition for $R_1$, we have
$\Ext^j_{C^\circ}(N\otimes_B P,R_1)=0$ for all $j<i$. 
Hence $\Ext^j_{B^\circ}(N,R)=0$ by Lemma \ref{yy4.3}(1).
This is the Auslander condition for $R$. The converse follows from 
the fact $R= R_1\otimes_C P^{\vee}$.

The argument above also shows that 
$$\opn{Cdim}_R(M)=\opn{Cdim}_{R_1}(M)\quad
\text{and}\quad \opn{Cdim}_{R}(N)=\opn{Cdim}_{R_1}(N\otimes_B P)$$
for all finite $A$-modules $M$ and finite $B^\circ$-modules $N$.
Since $\opn{GKdim}$ is preserved by Morita
equivalence, $R$ is Cohen-Macaulay if and only if $R_1$ is.
\end{proof}

\begin{proof}[Proof of Corollary \tup{\ref{yy0.3}}]
By \cite{YZ}, $A$ has a dualizing complex $R$ in all cases. 
Furthermore in cases (i-iv), $R$ is Auslander and Cohen-Macaulay.
By Theorem \ref{yy0.2}, $R\cong P[n]$ for some invertible
$A$-bimodule $P$ and some integer $n$, and $A$ is a Gorenstein 
ring.   In  cases (i-iv) the Auslander-Gorenstein and 
Cohen-Macaulay  properties of $A$ follow from Proposition 
\ref{yy4.4}. 
\end{proof}


\section{Weyl Algebras in Positive Characteristics}

In this section we study dualizing complexes and two-sided tilting 
complexes over the Weyl algebras $A_n$ when $\opn{char} k>0$. 

\begin{prop}
\label{yy5.1} Let $B$ be an Azumaya algebra over its center
$\mrm{Z}(B)$, and suppose $\opn{Spec} \mrm{Z}(B)$ is connected. Let 
$A$ be another ring and
$T \in \msf{D}(\cat{Mod}(A \otimes B^\circ))$ a two-sided tilting 
complex. Then $T \cong P[n]$ for some integer $n$ 
and some invertible $A$-$B$-bimodule $P$.
\end{prop}

\begin{proof} 
Use the proof of \cite[2.7]{Ye2}, noting that 
for a prime ideal $\mfrak{p} \subset \mrm{Z}(B)$, the localization
$B \otimes_{\mrm{Z}(B)} \mrm{Z}(B)_{\mfrak{p}}$ 
is a local ring.
\end{proof}

The following lemma takes care of dualizing complexes.

\begin{lem}
\label{yy5.2} Let $A$ be a left noetherian ring and $B$ a right
noetherian ring.
\begin{enumerate}
\item{} If $A$ has finite injective dimension as left module, and 
$B$ has finite injective dimension as right module, then
every two-sided tilting complex $T$ over $(A, B)$ is also a  
dualizing complex over $(A,B)$.
\item{} If $A$ or $B$ is noetherian and Gorenstein, then every 
dualizing complex over $(A,B)$ is also a two-sided tilting complex.
\end{enumerate}
\end{lem}

\begin{proof} (1) By \cite[1.6 and 1.7]{Ye2} 
the cohomologies of $T$ are finite modules on both sides and the 
morphisms $B\to \mrm{R}\opn{Hom}_A(T,T)$ and 
$A\to \mrm{R} \opn{Hom}_{B^\circ}(T,T)$ are isomorphisms.
Since $T$ has finite projective dimension over $A$ and the left 
module $A$ has finite injective dimension it follows that
$T$ also has finite injective dimension over $A$. 
Likewise on the right.

\medskip \noindent
(2) If $A$ is a noetherian Gorenstein ring then the bimodule $A$ 
is a dualizing complex over $A$. Let $R$ be any dualizing
complex over $(A,B)$. As mentioned earlier, the proof of 
\cite[4.5]{Ye2} -- suitably modified to fit our situation -- 
shows that the complex
$\mrm{R} \opn{Hom}_A(A,R)$ is a two-sided tilting complex over 
$(A, B)$. But $R \cong \mrm{R} \opn{Hom}_A(A,R)$.
\end{proof}

\begin{prop}
\label{yy5.3}
Let $B$ be an Azumaya algebra over its center $\mrm{Z}(B)$. 
Suppose $\mrm{Z}(B)$ is a noetherian Gorenstein ring and 
$\opn{Spec} \mrm{Z}(B)$ is connected.
Let $A$ be a left noetherian ring and $R$ a dualizing complex over 
$(A, B)$. Then $R \cong P[n]$ for some
integer $n$ and some invertible $A$-$B$-bimodule $P$.
\end{prop}

\begin{proof}
First we show that $B$ is also Gorenstein. Let $d$ be the 
injective dimension of $C := \mrm{Z}(B)$. For any prime ideal 
$\mfrak{p}$ of $C$ the local ring 
$C_{\mfrak{p}}$ is Gorenstein, of injective dimension $\leq d$,
and hence also the completion 
$\widehat{C}_{\mfrak{p}}$. Now the completion
$\widehat{B}_{\mfrak{p}} := B \otimes_{C} \widehat{C}_{\mfrak{p}}$
is isomorphic to a matrix ring 
$\mrm{M}_{r}(\widehat{C}_{\mfrak{p}})$, 
so by Morita equivalence $\widehat{B}_{\mfrak{p}}$ is Gorenstein. 
Faithful flatness (going over all primes $\mfrak{p}$) shows the 
vanishing of 
$\opn{Ext}^{i}_{B}(M, B)$ and $\opn{Ext}^{i}_{B^{\mrm{o}}}(N, B)$
for all finite modules $M$ and $N$ and all $i > d$, 
so we deduce that $B$ is Gorenstein.

Now we may use Proposition \ref{yy5.1} and Lemma \ref{yy5.2}(2).
\end{proof}

\begin{cor} 
Let $B$ be the $n$th Weyl algebra over $k$, with 
$\opn{char} k > 0$. Let $A$ be any left noetherian $k$-algebra, 
and let $R$ be any dualizing complex, or any two-sided tilting 
complex, over $(A, B)$. Then $R \cong P[n]$ for some invertible 
$A$-$B$-bimodule $P$ and integer $n$.
\end{cor}

\begin{proof}
By a result of Revoy \cite{Re}, the Weyl algebra $B$ is Azumaya 
with center a polynomial algebra over $k$. Now use the  
Propositions \ref{yy5.1} and \ref{yy5.3}.
\end{proof}



%
%
%
%
%
%


\section{Goodearl-Warfield's example}

We use an example of Goodearl-Warfield \cite[4.6]{GW2} to
show that not every noetherian simple ring has finite injective
dimension. This can also be done for the example of Shamsuddin
\cite{Sh}. 

\begin{exa}
\label{yy6.1}
Let $R[\theta;\delta]$ be the noetherian simple domain of infinite
Krull dimension constructed in \cite[4.6]{GW2}. In this example the 
base field $k$ is an infinite extension of $\mathbb Q$. 
The ring $R$ is a noetherian regular commutative $k$-algebra of 
infinite Krull dimension obtained by localizing a polynomial 
ring of countably many variables, which is essentially the example 
of Nagata \cite[Example 1, p.203]{Na}. 

Let $\fp$ be a prime ideal of $R$ of arbitrarily large height.
Then $R_{\fp}$ has finite global dimension, which is equal to
the height of $\fp$. Hence $R_{\fp}[\theta;\delta]$ has finite global
dimension, which is equal to the global dimension of $R_{\fp}$
\cite[2.1 and 3.4]{GW2}. Therefore $R_{\fp}[\theta;\delta]$ can have 
arbitrarily large injective dimension. Since $R_{\fp}[\theta;\delta]$ 
is a localization of $R[\theta;\delta]$ \cite[1.1]{GW2}, the injective 
dimension of $R[\theta;\delta]$ must be infinite.

By Theorem \ref{yy0.2} there is no dualizing complex over $R[\theta;\delta]$. 
\end{exa}

We conclude this paper by the following question.

\begin{que} 
\label{yy6.2}
Does every noetherian finitely generated simple ring of finite
Krull (or Gelfand-Kirillov) dimension have finite left and right 
injective dimension?
\end{que}

\bigskip
\centerline{\bf Acknowledgments}
\bigskip

Both authors were supported by the US-Israel Binational
Science Foundation. The second author was also supported by the 
National Science Foundation. Authors thank Ken Goodearl and
Paul Smith for several conversations on the subject.


\end{document}